\documentclass[12pt]{article}
\usepackage{mathrsfs}
\usepackage{epic,eepic,epsf,epsfig}
\usepackage{amsfonts,srcltx,mathrsfs}
\usepackage{multirow}
\usepackage{amsmath}
\usepackage{amssymb}
\usepackage{amsbsy}
\usepackage{graphicx}
\usepackage{amsfonts}
\usepackage{color}
\newtheorem{lem}{Lemma}[section]%
\newtheorem{theorem}[lem]{Theorem}%
\newtheorem{cor}[lem]{Corollary}%
\newtheorem{prop}[lem]{Proposition}%

\def\a{\alpha}

 \def\O{\Omega} \def\G{\Gamma}

\def\di{\bigm|}  
\def\nd{\mathrel{\bigm|\kern-.7em/}}

\def\f{\noindent}

\def\PSL{\hbox{\rm PSL}}\def\PSU{\hbox{\rm PSU}}
  \def\GL{\hbox{\rm GL}} \def\Mult{\hbox{\rm Mult}}
\def\PSp{\hbox{\rm PSp}}\def\P\GammaL{\hbox{\rm P\Gamma L}}

\def\Aut{\hbox{\rm Aut}}

\def\Cos{\hbox{\rm Cos}}

\newcommand{\qed}{\mbox{\raisebox{0.7ex}{\fbox{}}} \vspace{4truemm}}
\def\mz{{\mathbb Z}}

\begin{document}
\title{Heptavalent symmetric graphs with solvable stabilizers admitting vertex-transitive non-abelian simple groups}

\author{Jia-Li Du, Yan-Quan Feng\footnotemark, Yu-Qin Liu\\
{\small\em Department of Mathematics, Beijing
Jiaotong University, Beijing 100044, China}}

\footnotetext[1]{Corresponding author. E-mails:
JiaLiDu@bjtu.edu.cn, yqfeng@bjtu.edu.cn, 14121547@bjtu.edu.cn}

\date{}
 \maketitle

\begin{abstract}
A graph $\G$ is said to be {\em symmetric} if its automorphism group $\Aut(\G)$ acts transitively on the arc
set of $\G$. In this paper, we show that if $\G$ is a finite connected heptavalent symmetric graph
with solvable stabilizer admitting a vertex-transitive non-abelian simple group $G$ of automorphisms,
then either $G$ is normal in $\Aut(\G)$, or $\Aut(\G)$ contains a non-abelian simple normal subgroup $T$ such that $G\leq T$ and $(G,T)$ is explicitly given as one of $11$ possible exception pairs of non-abelian simple groups. Furthermore, if $G$ is regular on the vertex set of $\G$ then the exception pair $(G,T)$ is one of $7$ possible pairs, and if $G$ is arc-transitive then the exception pair $(G,T)=(A_{17},A_{18})$ or $(A_{35},A_{36})$.

\bigskip
\f {\bf Keywords:} Symmetric graph, Cayley graph, coset graph, simple group.

\medskip
\f {\bf 2010 Mathematics Subject Classification:} 05C25, 20B25.

\end{abstract}

\section{Introduction}
Throughout this paper, all groups and graphs are finite, and all graphs are simple and undirected.
Let $G$ be a permutation group on a set $\O$ and let $\a\in \O$. Denote by $G_{\a}$ the stabilizer of $\a$ in $G$, that is, the subgroup of $G$ fixing the point $\a$. The group $G$ is {\it semiregular} if $G_{\a}=1$ for every $\a\in\O$, and {\it regular} if $G$ is transitive and semiregular.
Denote by $\mz_n$, $D_n$, $F_n$ and $A_n$, the cyclic group of order $n$, the dihedral group of order $2n$,
the Frobenius group of order $n$ and the alternating group of degree $n$, respectively.

For a graph $\G$, we denote its vertex set and automorphism group by $V(\G)$ and $\Aut(\G)$, respectively.
The graph $\G$ is said to be {\it $G$-vertex-transitive} for $G\leq \Aut(\Gamma)$ if $G$ acts transitively on $V(\G)$, {\it $G$-regular} if $G$ acts regularly on $V(\G)$, and {\it $G$-arc-transitive} if $G$ acts transitively on the arc set of $\G$ (an arc is an ordered pair of adjacent vertices). In particular, $\G$ is {\it vertex-transitive} or {\it symmetric} if it is $\Aut(\G)$-vertex-transitive or $\Aut(\G)$-arc-transitive, respectively.

Let $G$ be a non-abelian simple group. The two extreme cases of $G$-vertex-transitive graph are: $G$-arc-transitive graph and $G$-regular graph. Let $\G$ be a connected cubic $G$-arc-transitive graph. Li~\cite{CHLi} proved that either $G\unlhd \Aut(\G)$ or $(G,\Aut(\G))=(A_7,A_8)$, $(A_7,S_8)$, $(A_7,2.A_8)$, $(A_{15},A_{16})$ or $(\GL(4,2),\rm AGL(4,2))$. Fang {\em et al}~\cite{FangLW} proved that none of the above five pairs can happen, that is, $G$ is always normal in $\Aut(\G)$. Du {\em et al}~\cite{DFZh} showed that a connected pentavalent $G$-arc-transitive graph $\G$ is either $G\unlhd \Aut(\G)$ or $\Aut(\G)$ contains a non-abelian simple normal subgroup $T$ such that $G\leq T$ and $(G,T)$ is one of $17$ possible pairs of non-abelian simple groups.

A $G$-regular graph is also called a {\em Cayley graph} over $G$, and the Cayley graph is called {\em normal} if $G\unlhd \Aut(\G)$. One of the most remarkable achievements about Cayley graphs over $G$ is the complete classification of connected cubic symmetric non-normal Cayley graphs over $G$. In 1996, Li~\cite{CHLi} proved that a connected cubic symmetric Cayley graph $\G$ over $G$ is either normal or $G=A_5$, $A_7$, $\PSL(2,11)$, $M_{11}$, $A_{11}$, $A_{15}$, $M_{23}$, $A_{23}$ or $A_{47}$. In 2005, Xu {\em et al}~\cite{XFWX2005} proved that either $\G$ is normal or $G= A_{47}$, and two years later, Xu {\em et al}~~\cite{XFWX} further showed that if $G=A_{47}$ and $\G$ is not normal, then $\G$ must be $5$-arc-transitive and up to isomorphism there are exactly two such graphs. The result is also proved by Conder~\cite{Conder} and Li~{\em et al}~\cite{LL}. In 2005, Du {\em et al} \cite{DW} gave a classification of the connected symmetric cubic Cayley graphs on $\PSL(2,p)$ with $p\geq 5$ a prime. Du {\em et al}~\cite{DFZh} showed that a connected pentavalent symmetric Cayley graph $\G$ over $G$ is either normal, or $\Aut(\G)$ contains a non-abelian simple normal subgroup $T$ such that $G\leq T$ and $(G,T)$ is one of $13$ possible pairs of non-abelian simple groups.

In this paper, we investigate connected heptavalent symmetric graphs with solvable stabilizers admitting
vertex-transitive non-abelian simple groups.

\begin{theorem}\label{theo=main}
Let $G$ be a non-abelian simple group and let $\G$ be a connected heptavalent symmetric $G$-vertex-transitive graph
with solvable stabilizer. Then either $G\unlhd \Aut(\G)$, or $\Aut(\G)$ contains a non-abelian simple normal subgroup $T$ such that $G\leq T$ and
$(G,T)=(A_6,A_7)$ or $(A_{n-1},A_n)$ with $n\geq 14$ and $n\di 2^2\cdot3^2\cdot7$.
\end{theorem}

\begin{cor}\label{cor=arc}
Let $G$ be a non-abelian simple group and let $\G$ be a connected heptavalent $G$-arc-transitive graph with solvable stabilizer. Then either $G\unlhd \Aut(\G)$, or $\Aut(\G)$ contains a non-abelian simple normal subgroup $T$ such that $G\leq T$ and $(G,T)=(A_{17},A_{18})$ or $(A_{35},A_{36})$.
\end{cor}

\begin{cor}\label{cor=regular}
Let $G$ be a non-abelian simple group and let $\G$ be a connected heptavalent symmetric $G$-regular graph with solvable stabilizer. Then either $G\unlhd \Aut(\G)$, or $\Aut(\G)$ contains a non-abelian simple normal subgroup $T$ such that $G\leq T$ and $(G,T)=(A_{n-1},A_n)$ with $n=7$, $3\cdot 7$, $2^2\cdot7$, $3^2\cdot7$, $2^2\cdot3\cdot7$, $2\cdot3^2\cdot7$ or $2^2\cdot3^2\cdot7$.
\end{cor}

\section{Preliminaries\label{s2}}

In this section, we describe some preliminary results which will be used later. First we describe solvable
vertex stabilizers of connected heptavalent symmetric graphs.

\begin{prop}\label{prop=stabilizer} {\rm \cite[Corollary 2.2]{Guo}}
Let $\Gamma$ be a connected heptavalent $G$-arc-transitive graph with $v\in V(\Gamma)$. If $G_v$ is solvable,\
then $G_v\cong \mathbb{Z}_7$, $D_{7}$, $F_{21}$, $D_{7}\times \mz_2$, $F_{21}\times \mz_3$, $F_{42}$,
$ F_{42}\times \mathbb{Z}_2$, $ F_{42}\times \mathbb{Z}_3$ or $ F_{42}\times \mathbb{Z}_6$, where $F_{n}$ is the Frobenius group of order $n$. In particular, $|G_v|\di 2^2\cdot3^2 \cdot 7$.
\end{prop}

The following result follows from the classification of three-factor simple groups.

\begin{prop}\label{prop=235simplegroup}{\rm \cite[Theorem \MakeUppercase{\romannumeral1} ]{Huppert2}}
Let $G$ be a non-abelian simple $\{2,3,7\}$-group.
Then $G= \PSL(2,7)$, $\PSL(2,8)$ or $\rm PSU(3,3)$.
\end{prop}

Let $G$ and $E$ be two groups. We call an extension $E$ of $G$ by $N$ a {\em central extension} of $G$ if $E$ has a central subgroup $N$ such that $E/N\cong G$, and if further $E$ is perfect, that is, the derived group $E'$ equals to $E$, we call $E$ a {\em covering group} of $G$. If the center has order $n$, the coving group is often referred to as a {\em n-cover}. Schur~\cite{Schur} proved that for every non-abelian simple group $G$ there is a unique maximal covering group $M$ such that every covering group of $G$ is a factor group of $M$ (see \cite[Kapitel V, \S23]{Huppert}). This group $M$ is called the {\em full covering group} of $G$, and the center of $M$ is the {\em Schur multiplier} of $G$, denoted by $\Mult(G)$.

By Kleidman and Liebeck~\cite{Kleidman}, $\Mult(A_n)=\mz_2$ for $n\geq 5$ with $n\not=6,7$, and $\Mult(A_n)=
\mz_6$ for $n= 6$ or $7$. This implies that $A_n$ has a unique $2$-cover for $n\geq 5$, denoted by $2.A_n$, and $A_7$ has a unique $3$-cover and $6$-cover, denoted by $3.A_7$ and $6.A_7$, respectively. It is well-known that for $n\geq 7$, $A_n$ has no proper subgroup of index less than $n$ and has exactly one conjugate class of subgroups of index $n$. This implies the following proposition.

\begin{prop}\label{prop=covering group} For $n\geq 7$, the alternating group $A_n$ has a unique $2$-cover $2.A_n$, and all subgroups of index $n$ of $2.A_n$ are conjugate and isomorphic to $2.A_{n-1}$.
Moreover, $A_7$ has a unique $3$-cover $3.A_7$ and a unique $6$-cover $6.A_7$, and all subgroups of index $7$ of $3.A_7$ and $6.A_7$ are conjugate and isomorphic to $3.A_6$ and $6.A_6$, respectively.
\end{prop}

Let $\G$ be a graph and $N \leq \Aut(\G)$. The \emph{quotient graph} $\G_N$ of $\G$ relative to $N$ is defined as the graph with vertices the orbits of $N$ on $V(\G)$ and with two orbits adjacent if there is an edge in $\G$ between these two orbits.

\begin{prop}\label{prop=atlesst3orbits}{\rm \cite[Theorem 9]{Lorimer}}
Let $\Gamma$ be a connected $G$-arc-transitive graph of prime valency, and let $N\unlhd G$ have at least three orbits on $V(\Gamma)$. Then $N$ is the kernel of $G$ on $V(\G_N)$, and semiregular on $V(\G)$. Furthermore, $\G_N$
is $G/N$-symmetric with $G/N \leq \Aut(\G_N)$.
\end{prop}

Let $G$ be a group. For $H\leq G$, let $D$ be a union of some double cosets of $H$ in $G$ such that $D^{-1}=D$. The \emph{coset graph} $ \Gamma=\Cos(G,H,D)$ on $G$ with respect to $H$ and $D$ is defined to have vertex set $V(\Gamma)=[G:H]$, the set of right cosets of $H$ in $G$, and edge set $E(\Gamma)=\{\{Hg,Hdg\}\ |\ g\in G,d\in D\}$. The graph $\Gamma$ has valency $|D|/|H|$ and it is connected if and only if $G=\langle D,H\rangle$, that is, $D$ and $H$ generate $G$. The action of $G$ on $V(\G)$ by right multiplication induces a transitive group of automorphisms, and this group is symmetric if and only if $D$ is a single double coset. Moreover, this action is faithful if and only if $H_G=1$, where $H_G$ is the largest normal subgroup of $G$ contained in $H$.

Conversely, let $\Gamma$ be a $G$-vertex-transitive graph with $G\leq \Aut(\Gamma)$. By~\cite{Sabidussi}, the graph $\Gamma$ is isomorphic to a coset graph $\Cos(G,H,D)$, where $H=G_{v}$ is the vertex stabilizer of $v\in V(\Gamma)$ in $G$ and $D$ consists of all elements of $G$ mapping $v$ to one of its neighbors.
It is easy to show that $H_G=1$ and $D$ is a union of some double cosets of $H$ in $G$ satisfying $D^{-1}=D$. Assume that $\G$ is $G$-arc-transitive and $g \in G$ interchanges $v$ and one of its neighbors. Then $g^2 \in H$ and $D=HgH$. Furthermore, $g$ can be chosen as a $2$-element in $G$, and the valency of $\Gamma$ is $|D|/|H|=|H:H\cap H^g|$. For more details regarding coset graph, referee to  \cite{FangP,Lorimer,Miller,Sabidussi}.

\begin{prop}\label{prop=cosetgraph}
Let $\Gamma$ be a connected $G$-arc-transitive graph of valency $k$, and let $\{u,v\}$ be an edge of $\Gamma$. Then $ \Gamma $ is isomorphic to a coset graph $\Cos(G,G_{v},G_{v}gG_{v})$, where $g$ is a $2$-element in $G$ such that $G_{uv}^{g}=G_{uv}$, $g^2\in G_{v}$,  $\langle G_v,g\rangle=G$, and $k=|G_v:G_v \cap G_{v}^{g}|$.
\end{prop}

In Proposition \ref{prop=cosetgraph}, the $2$-element $g$ such that $g^2\in G_{v}$,  $\langle G_v,g\rangle=G$ and $k=|G_v:G_v \cap G_{v}^{g}|$, is called {\em feasible} to $G$ and $G_v$. Feasible $g$ can be computed by MAGMA~\cite{magma} when the order $|G|$ and the valency $k$ are not too large, and this is used in Section~3.

\section{Proof of Theorem~\ref{theo=main}}

In this section, we always assume that $G$ is a non-abelian simple group and $\G$ is a connected heptavalent
symmetric $G$-vertex-transitive graph. Let $A=\Aut(\G)$ and $v\in V(\G)$. We further assume that $A_v$ is solvable.

\begin{lem}\label{lem=GH}
Let $H\leq A$ and $GH\leq A$. Then $|(GH)_v|/|G_v|=|H|/|H\cap G|$ is a divisor of $2^2\cdot 3^2\cdot7$.
\end{lem}

\f {\bf Proof:} Since $\G$ is $G$-vertex-transitive, by Frattini argument, $GH=G(GH)_v$ and
$|H||G|/$ $|H \cap G|=$ $|HG|=|G(GH)_v|=$ $|G||(GH)_v|/|G\cap (GH)_v|=$ $|G||(GH)_v|/|G_v|$.
Then $|H|/|H\cap G|=|(GH)_v|/|G_v|$. Since $GH\leq A$ and $A_v$ is solvable,
Proposition~\ref{prop=stabilizer} implies $|(GH)_v|\di 2^2\cdot 3^2\cdot7$.
It follows $|H|/|H\cap G|\di 2^2\cdot 3^2\cdot 7$.
\hfill\qed

The {\em radical} of a finite group is the largest solvable normal subgroup of the group.

\begin{lem}\label{lem=insolvable2}
Let $\G$ be $X$-arc-transitive with $X$ being of trivial radical.
Then either $G\unlhd X$ or $X$ has a non-abelian simple normal subgroup $T$ such that $G\leq T$
and $(G,T)=(A_6,A_7)$ or $(A_{n-1},A_n)$ with $n\geq 14$ and $n\di 2^2\cdot3^2\cdot7$.
\end{lem}

\f {\bf Proof:} Let $N$ be a minimal normal subgroup of $X$. Since $X$ has trivial radical, $N=T^s$ for a positive integer $s$ and a non-abelian simple group $T$. Then $NG\leq X$, and by Lemma~\ref{lem=GH}, $|N|/|N\cap G|\di 2^2\cdot 3^2\cdot 7$.

Suppose $N \cap G =1$. Then $|N|=|N|/|N\cap G|\di 2^2\cdot 3^2\cdot 7$. Write $K=NG$. Then $|K|=|NG|=|N||G|$ and $K/N\cong G$. Since $G$ is vertex-transitive, the Frattini argument implies $K=GK_v$, and so $|K|=|GK_v|=|G||K_v|/|G_v|$. It follows $|K_v|=|G_v||N|$. Since  $|\PSU(3,3)|=2^5 \cdot3^2\cdot 7$, Proposition~\ref{prop=235simplegroup} implies $N=\PSL(2,7)$ or $\PSL(2,8)$. Thus, $2^3\di |K_v|$, forcing that $2^3\di|X_v|$, which is impossible by Proposition \ref{prop=stabilizer}.

Thus, $N\cap G\not=1$. The simplicity of $G$ implies that $G\cap N=G$ and so $G\leq N$. Thus, $|N|/|G|=|N|/|N\cap G|\di 2^2\cdot 3^2\cdot 7$. If $G=N$, then $G\unlhd X$ and we have done. Now, we assume $G\neq N$, that is, $G$ is a proper subgroup of $N$, yielding $N_v\not=1$. Since $X$ is arc-transitive, $X_v$ is primitive on the $7$ neighbours of $v$ in $\G$, and since $N_v\unlhd X_v$, we have $7\mid |N_v|$. In particular, $\G$ is $N$-symmetric. Recall that $N=T^s$ for some $s\geq 1$.

Suppose $s\geq 2$. Since $T\unlhd N$, we have $T\cap G=1$ or $G$. If $T\cap G=G$ then $G\leq T$ and so $|G|\di |T|$. Thus, $|T|^{s-1}|T:G|\di |N|/|G|$, implying $|T|\di |N|/|G|$. If $T\cap G=1$, we also have $|T|\di |N|/|G|$ as $|T|=|T|/|G\cap T|=|TG|/|G|\di |N|/|G|$. By Lemma~\ref{lem=GH}, $|T|\di 2^2\cdot3^2\cdot7$. Since $G\leq N$, both $G$ and $T$ are non-abelian simple $\{2,3,7\}$-groups. Proposition~\ref{prop=235simplegroup} implies $G=\PSL(2,7)$ or $\PSL(2,8)$, and also $T=\PSL(2,7)$ or $\PSL(2,8)$. Since $7^2\nmid |N|/|G|$, we have $s=2$ and $|G|\di |T|$. If both direct factors $T$ of $N$ are not semiregular, then $7\di |T_v|$ for each $T$ because $\G$ is $N$-symmetric, yielding $7^2\di |N_v|$, contrary to Proposition~\ref{prop=stabilizer}. Thus, $N$ has a direct factor $T$ with $T_v=1$. Since $|G|\di |T|$ and $\G$ is $G$-vertex-transitive, $T$ is regular on $V(\Gamma)$ and so $|N_v|=|T|=|\PSL(2,7)|$ or $|\PSL(2,8)|$, implying $2^3\di |N_v|$, which is impossible by Proposition~\ref{prop=stabilizer}.

The above contradiction implies $s=1$, that is, $N$ is the non-abelian simple group $T$. Thus, $G\leq T\unlhd X$
and $|T:G|=|N|/|G|\di 2^2\cdot 3^2\cdot 7$. If $G=T$ then $G\unlhd X$ and we are done. Assume $G<T$, that is, $G$ is a proper subgroup of $T$. Then $T_v\not=1$ and so $7\di |T_v|$. In particular, $7\di |T|$.
Let $M$ be a maximal subgroup of $T$ containing $G$. Then $|T:G|=|T:M||M:G|$, and $T$ can be viewed as a primitive permutation group of degree dividing $2^2\cdot 3^2\cdot7$ with a stabilizer containing the non-abelian simple group $G$. By \cite[Section 7]{Colva}, $T$ and $M$ are listed in Table~\ref{table=2} (note that $M$ is not solvable):

\begin{table}[ht]
\begin{center}
\begin{tabular}{|c|c|c|c|c|c|}

\hline
$T$       & $M$       & $|T:M|$              & $T$          & $M$             & $|T:M|$      \\
\hline
 $A_7$    & $S_5$     & $3\cdot7$           &$A_8$          & $S_6$           & $2^2\cdot7$   \\
\hline
$A_9$          & $S_7$    & $2^2\cdot3^2$    &$A_9$   & $(A_4\times A_5)\rtimes \mz_2$   & $2\cdot 3^2\cdot7$ \\
\hline
$A_9$          & $(A_6\times \mz_3)\rtimes \mz_2$    & $2^2\cdot3\cdot7$  & $\PSU(3,3)$   & $\PSL(2,7)$        & $2^2\cdot 3^2$  \\
\hline
$\PSU(4,2)$  & $S_6$   & $2^2\cdot 3^2$       & $\PSU(4,3)$      & $\PSU(4,2)$      & $2\cdot3^2\cdot7$ \\
\hline
$\PSp(6,2)$  & $S_8$      & $2^2\cdot 3^2$    & $\PSp(6,2)$      & $\PSU(4,2)\rtimes \mz_2$    & $2^2\cdot 7$ \\
\hline
$\PSp(6,2)$   & $\mz_2^5\rtimes S_6$  & $ 3^2\cdot 7$     & $\PSL(3,4)$   & $\mz_2^4\rtimes A_5$  &$3\cdot 7$ \\
\hline
$M_{12}$   & $M_{11}$  & $ 2^2\cdot 3$     & $M_{11}$   & $\PSL(2,11)$  &$2^2\cdot 3$ \\
\hline
$A_n$       & $A_{n-1}$   & $n\di 2^2\cdot3^2\cdot7$           &    &  & \\
\hline
\end{tabular}
\end{center}
\vskip -0.5cm
\caption{{Primitive group of non-abelain simple groups of degree dividing $2^2\cdot3^3\cdot7$}}\label{table=2}
\end{table}

Since $7\di |T|$, we have $(T,M)\not=(M_{12},M_{11}),(M_{11},\PSL(2,11)), (A_6,A_5)$.

Suppose $(T,M)=(A_8,S_6)$. Then $|M:G|\di 3^2\cdot7$. Since $M$ has no subgroup of index less than $6$
and $|M|=2^4\cdot 3^2\cdot5$, we have $|M:G|=9$, which is impossible because $S_6$ has no subgroup of index $9$. Thus, $(T,M)\not=(A_8,S_6)$, and similarly, since $T$ has a non-abelian simple subgroup $G$ with $|T:G|\di 2^2\cdot 3^2\cdot7$, we have
$(T,M)\not=(A_9$, $S_7)$, $(A_9$, $(A_4\times A_5)\rtimes \mz_2)$, $(A_9$, $(A_6\times \mz_3)\rtimes \mz_2)$, $(\PSp(6,2)$, $(\PSU(4,2)\rtimes \mz_2))$, $(\PSp(6,2)$, $S_8)$, $(\PSp(6,2)$, $\mz_2^5\rtimes S_6)$, $(\PSL(3,4)$, $\mz_2^5\rtimes A_5)$ or $(\PSU(4,2)$, $S_6)$.

Let $(T,M)=(A_7,S_5)$. Then $|M:G|\di 2^2\cdot3$. Since $G$ is non-abelian simple, we have $G=A_5$, and hence $(T,G)=(A_7,A_5)$. Let $(T,M)=(\PSU(3,3),\PSL(2,7))$ or $(\PSU(4,3),\PSU(4,2))$.
Then $G$ is a non-abelain simple subgroup of $M$ with index dividing $7$ or $2$, respectively. It follows $M=G$ and so $(T,G)=(\PSU(3,3),\PSL(2,7))$, $(\PSU(4,3),\PSU(4,2))$.

Let $(T,M)=(A_n,A_{n-1})$ with $n\di 2^2\cdot 3^2\cdot7$ and $n\geq 7$.
If $G=M$, then $(T,G)=(A_n,A_{n-1})$. If $G<M$, then $|M:G|\geq n-1$ and so $|T:G|\geq n(n-1)$.
It follows $n<16$ as $17\cdot16>252$. Since $n\di 2^2\cdot 3^2\cdot7$, we have $n=7,9,12$ or $14$, and by Atlas~\cite{Atlas}, the only possible pair is $(T,M,G)=(A_7,A_6,A_5)$ and thus $(T,G)=(A_7,A_5)$.

Now we obtain $(T,G)=(\PSU(3,3)$, $\PSL(2,7))$, $(\PSU(4,3)$, $\PSU(4,2))$, $(A_7$, $A_5)$ or $(A_n$, $A_{n-1})$ with $n\geq 7$ and $n\di 2^2\cdot 3^2\cdot7$. To finish the proof, we only need to show $(T,G)\neq $ $(\PSU(3,3)$, $\PSL(2,7))$, $(\PSU(4,3)$, $\PSU(4,2))$, $(A_7$, $A_5)$, $(A_9$, $A_8)$ or $(A_{12}$, $A_{11})$.

Suppose $(T,G)=(\PSU(3,3),\PSL(2,7))$. Then $|T_v|/|G_v|=|T|/|G|= 2^2\cdot3^2$.
By Proposition \ref{prop=stabilizer}, $|T_v|=2^2\cdot 3^2\cdot7$
and so $|T:T_v|=2^3\cdot3$. However, by Atlas \cite[pp. 14]{Atlas},
$\PSU(3,3)$ has no subgroup with index $24$, a contradiction.

Suppose $(T,G)=(\PSU(4,3),\PSU(4,2))$. Then $|T_v|/|G_v|=|T|/|G|= 2\cdot3^2\cdot 7$.
By Proposition~\ref{prop=stabilizer}, $|T_v|=2\cdot3^2\cdot 7$ or $2^2\cdot 3^2\cdot7$, and $T_v\cong F_{42}\times \mz_3$ or $F_{42}\times \mz_6$. However, by MAGMA, $\PSU(4,3)$ has no such subgroups, a contradiction.

Suppose $(T,G)=(A_7,A_5)$. Then $|T_v|/|G_v|=|T|/|G|=2\cdot3\cdot 7$ and $2\cdot 3\cdot 7\di|T_v|$. By Proposition~\ref{prop=stabilizer}, $T_v$ contains a subgroup isomorphic to $F_{42}$, but this is impossible because   $A_7$ has no such a subgroup by Atlas \cite[pp. 10]{Atlas}.

Suppose $(T,G)=(A_9,A_8)$. Then $|T_v|/|G_v|=|T|/|G|=3^2$ and $3^2\di|T_v|$.
By Proposition \ref{prop=stabilizer}, $T_v\cong F_{21}\times \mz_3$, $F_{42}\times \mz_3$ or
$F_{42}\times \mz_3$. However, by MAGMA, $A_9$ has no such subgroups, a contradiction.

Suppose $(T,G)=(A_{12},A_{11})$. Then $|T_v|/|G_v|=|T|/|G|=2^2\cdot3$.
By Proposition \ref{prop=stabilizer}, $T_v\cong F_{42}\times \mz_2$ or
$F_{42}\times \mz_6$. By MAGMA, $T$ has no subgroup isomorphic to
$F_{21}\times \mz_6$ and has one conjugacy class of subgroups isomorphic to $F_{42}\times \mz_2$. Take a given $T_v$ in the class. By Proposition \ref{prop=cosetgraph}, there is a $2$-element $g$ in $T$ such that $g^2\in T_{v}$,  $\langle T_v,g\rangle=T$ and $7=|T_v:T_v \cap T_{v}^{g}|$, but computation shows that there is no such $g$ in $T$, a contradiction.
\hfill\qed

\begin{lem}\label{lem=T normal}
Let $A$ have a non-trivial radical $R$ with at least three orbits.
Assume $G\ntrianglelefteq A$ and $RG=R\times G$. Then $A$ contains a non-abelian simple
normal subgroup $T$ such that $G\leq T$ and $(G,T)=(A_6,A_7)$ or $(A_{n-1},A_n)$ with $n\geq 14$ and $n\di 2^2\cdot3^2\cdot7$.
\end{lem}

\f {\bf Proof:} Set $B=RG=R\times G$. Since $R\neq 1$ has at least three orbits,
by Proposition~\ref{prop=atlesst3orbits}, the quotient graph $\G_{R}$ is a connected heptavalent
$A/R$-arc-transitive graph with $A/R\leq \Aut(\G_{R})$.
Furthermore, $\G_R$ is $B/R$-vertex-transitive and $G$
is characteristic in $B$. It follows that $B\ntrianglelefteq A$ because $G\ntrianglelefteq A$, and
hence $G\cong B/R\ntrianglelefteq A/R$. Since $R$ is the largest solvable normal subgroup of $A$,
$A/R$ has trivial radical. Let $\a$ be the orbit of $R$ containing $v$. Then $A_\a=RA_v$ and $(A/R)_\a=A_vR/R$. It follow that $(A/R)_\a$ is solvable as $A_v$ is solvable. By Lemma~\ref{lem=insolvable2}, $A/R$ has a non-abelian simple normal subgroup $I/R$ such that $B/R\leq I/R\cong T$ and $(B/R,I/R)\cong (G,T)=(A_6,A_7)$ or $(A_{n-1},A_n)$ with $n\geq 14$ and $n\di 2^2\cdot3^2\cdot7$.

Let $C=C_I(R)$. Then $G\leq C$ as $B=R\times G\leq I$. Clearly, $C\unlhd I$ and $C\cap R\leq Z(C)$.
Since $G\cong (R\times G)/R\leq CR/R$ and $C/C\cap R\cong CR/R\unlhd I/R\cong T$, we have  $C\cap R=Z(C)$, $C/Z(C)\cong T$
and $I=CR$. Thus, $C'/C'\cap Z(C)\cong C'Z(C)/Z(C)=(C/Z(C))'=C/Z(C)\cong  T$, and so  $Z(C')=C'\cap Z(C)$, $C=C'Z(C)$ and $C'/Z(C')\cong T$. Furthermore, $C'=(C'Z(C))'=C''$, implying that $C'$ is a covering group of $T$. It follows $G\leq C'$ because $C/C'$ is abelian.

Recall that $T=A_7$ or $A_n$ with  $n\geq 14$ and $n\di 2^2\cdot3^2\cdot7$. By \cite[Theorem 5.1.4]{Kleidman} and the Atlas~\cite[p. 123]{Atlas}, $\Mult(A_7)=\mz_6$ and $\Mult(A_n)=\mz_2$ for $n\geq 8$. Since $C'$ is a covering group of $T$, we have $Z(C')=1$, $\mz_2$, $\mz_3$ or $\mz_6$.

Suppose $Z(C')=\mz_2$. Then $G\times Z(C')$ is a subgroup of index $n$ of $C'$ isomorphic to $A_{n-1}\times \mz_2$, which is impossible by Proposition~\ref{prop=covering group}.

Suppose $Z(C')=\mz_3$ or $\mz_6$. Then $(G,T)=(A_6,A_7)$ and $G\times Z(C')$ is a subgroup of index $7$ of $C'$ isomorphic to $A_{6}\times \mz_3$ or $A_6\times \mz_6$, which are impossible by Proposition~\ref{prop=covering group}.

Thus, $Z(C')=1$. Furthermore, $C'\cong T$ and $G<C'\unlhd I$. Since $|I|=|I/R||R|=|T||R|=|C'||R|$ and $C'\cap R=1$, we have $I=C'\times R$. This implies that $C'$ is characteristic in $I$, and so $C'\unlhd A$ because $I\unlhd A$. It follows that $A$ has a non-abelian simple normal subgroup $C'$ such that $G\leq C'$ and $(G,C')\cong (G,T)=(A_6,A_7)$ or $(A_{n-1},A_n)$ with $n\geq 14$ and $n\di 2^2\cdot3^2\cdot7$. This completes the proof.
\hfill\qed

Now, we are ready to prove Theorem \ref{theo=main}.

\medskip

\f {\bf The proof of Theorem \ref{theo=main}:} Let $G$ be a non-abelian simple group and let $\G$ be a connected heptavalent symmetric $G$-vertex-transitive graph. Write $A=\Aut(\G)$ and let $R$ be the radical of $A$. To prove the theorem, we may assume $G\ntrianglelefteq A$.
If $R=1$, the theorem is true by  Lemma~\ref{lem=insolvable2}. Now we assume $R\not=1$.

Set $B=RG$. Then $G\cap R=1$ and $|B|=|R||G|$. By Lemma~\ref{lem=GH},
$|R|=|B|/|G|=|B_v|/|G_v|\di 2^2\cdot 3^2\cdot 7$.

If $R$ is transitive on $V(\G)$, Lemma \ref{lem=GH} implies $|G|=|B|/|R|=|(B)_v|/|R_v|\di 2^2\cdot 3^2\cdot 7$, which is impossible by Proposition \ref{prop=235simplegroup}. If $R$ has two orbits, then $\G$ is bipartite and so $G$ has a subgroup of index $2$, which is impossible. Thus, $R$ has at least three orbits.

To finish the proof, by Lemma~\ref{lem=T normal}, it suffices to show $B=R\times G$.

By Proposition~\ref{prop=atlesst3orbits}, $R$ is semiregular on $V(\G)$ and the quotient graph $\G_{R}$ is a
connected heptavalent $A/R$-symmetric graph with $A/R\leq \Aut(\G_{R})$.  Moreover, $\G_R$ is $B/R$-vertex-transitive. Since $R$ is the largest solvable normal subgroup of $A$, $A/R$ has trivial radical.
Since $A_v$ is solvable, $(A/R)_{\a}$ is solvable with $\a\in V(\G_R)$.
By Lemma~\ref{lem=insolvable2}, $B/R\unlhd A/R$ or $A/R$ has a non-abelian simple normal subgroup $I/R$ such that $B/R\leq I/R$ and $(B/R,I/R) \cong (G,T)$ with $(G,T)=(A_6,A_7)$ or $(A_{n-1},A_n)$ with $n\geq 14$ and $n\di 2^2\cdot3^2\cdot7$.

Since $|R|\di 2^2\cdot 3^2\cdot 7$, we may write $|R|=2^m\cdot3^n\cdot 7^k$, where $0\leq m\leq 2$,
$0\leq n\leq 2$ and $0\leq k\leq 1$. Since $R$ is solvable, there exists a series of subgroups of $B$:
\begin{center}
$B>R=R_{s}>\cdots >R_1>R_0=1$
\end{center}
such that $R_i\unlhd B$ and $R_{i+1}/R_i$ is an elementary abelian $r$-group with $0\leq i\leq s-1$, where
$r=2$, $3$ or $7$. Clearly, $G\leq B$ has a natural action on $R_{i+1}/R_i$ by conjugation.

Suppose $B\neq R\times G$. Then there exists some $0\leq j\leq s-1$ such that $GR_i=G\times R_i$ for any $0\leq i\leq j$,
but $GR_{j+1}\not=G\times R_{j+1}$. If $G$ acts trivially on $R_{j+1}/R_j$ by conjugation, then $[GR_j/R_j,R_{j+1}/R_j]=1$.
Since $GR_j/R_j\cong G$ is simple, we have $(GR_j/R_j)\cap (R_{j+1}/R_j)=1$. Note that $|GR_{j+1}/R_j|=|GR_{j+1}/R_{j+1}||R_{j+1}/R_j|=|G||R_{j+1}/R_j|=|GR_j/R_j||R_{j+1}/R_j|$. Then $GR_{j+1}/R_j=GR_j/R_j\times R_{j+1}/R_j$. In particular, $GR_j\unlhd GR_{j+1}$ and so $G\unlhd GR_{j+1}$ because $GR_j=G\times R_j$
implies that $G$ is characteristic in $GR_j$. It follows that $GR_{j+1}=G\times R_{j+1}$, a contradiction.
Thus, $G$ acts non-trivially on $R_{j+1}/R_j$, and the simplicity of $G$ implies that $G$ acts faithfully
on $R_{j+1}/R_j$.

Set $|R_{j+1}/R_j|=r^\ell$. Since $|R|=2^m\cdot3^n\cdot 7^k$, we have $\ell\leq m\leq 2$ for $r=2$ or $3$
and $\ell\leq 1$ for $r=7$. Recall that $G$ acts faithfully on $R_{j+1}/R_j$ and $R_{j+1}/R_j$ is elementary abelian. Then $G\leq \GL(\ell,r)$. For each case, $\GL(\ell,r)$ is solvable and so $G$ is solvable, a contradiction.

The above contradiction implies $B=R\times G$, as required.
\hfill\qed

\f {\bf The proof of Corollary \ref{cor=arc}:} Let $G$ be a non-abelian simple group and let $\G$ be a connected heptavalent $G$-arc-transitive graph with $\Aut(\G)_v$ being of solvable. Assume $G\ntrianglelefteq \Aut(\G)$. By Theorem~\ref{theo=main},
$\Aut(\G)$ has a non-abelian simple normal subgroup $T$ such that $G\leq T$ and $(G,T)=(A_6,A_7)$ or $(A_{n-1},A_n)$ with $n\geq 14$ and $n\di 2^2\cdot3^2\cdot7$.

Since $G$ is symmetric, both $|G_v|$ and $ |T_v|$ are divisible by $7$, and so $|T|/|G|=|T_v|/|G_v|$ is a divisor of $2^2\cdot3^2$ by Proposition~\ref{prop=stabilizer}. Thus, $(G,T)=(A_{n-1},A_n)$ with $n=6,9,12,18$ or $36$, and by Theorem~\ref{theo=main}, $(G,T)=(A_{n-1},A_n)$ with $n=18$ or $36$. \hfill\qed

\f {\bf The proof of Corollary \ref{cor=regular}:} Let $G$ be a non-abelian simple group and let $\G$ be a connected heptavalent $G$-regular graph with $\Aut(\G)_v$ being of solvable. Assume $G\ntrianglelefteq \Aut(\G)$. By Theorem~\ref{theo=main},  $\Aut(\G)$ has a non-abelian simple normal subgroup $T$ such that $G\leq T$ and $(G,T)$ is as given in Theorem~\ref{theo=main}. Since $G$ is regular and $G\leq T\unlhd A$, we have $|T_v|=|T|/|G|$, and by Proposition~\ref{prop=stabilizer}, $(G,T)=(A_{n-1},A_n)$ with $n=7$, $2\cdot 7$, $3\cdot7$, $2^2\cdot7$, $3^2\cdot 7$, $2\cdot 3\cdot7$, $2^2\cdot 3\cdot7$, $2\cdot 3^2\cdot7$ or $2^2\cdot 3^2\cdot7$.
To finish the proof, it suffices to show that $n\neq 2\cdot7, 2\cdot3\cdot7$. Now suppose to the contrary that $n=2\cdot7$ or $2\cdot3\cdot7$.

By Proposition~\ref{prop=cosetgraph}, $\G=\Cos(T,H,HgH)$ with $H=T_v$ for some $2$-element $g$. Since $G$ is regular, we have $T=GH$ with $
G\cap H=1$. Then $H$ acts regularly on $\Omega=[T:G]=\{Gt|t\in T\}$ by right multiplication. Note that $|H|=|\Omega|=n$.
Write $\Omega=\{\delta_1, \delta_2,\cdots ,\delta_n \}$ with $\delta_1=G$ and set $L=H\cap H^g$. Since $H$ is transitive on the neighbours of $v$ in $\G$, we have $|H:L|=7$, and the regularity of $H$ on $\Omega$ implies that $L$ has exactly seven orbits on $\Omega$, say $\Omega_i=\{\delta_{ij}\ |\ j=1,2,\cdots, n/7\}$, $1\leq i\leq 7$.

By Proposition~\ref{prop=stabilizer}, $H=D_7$ or $F_{42}$ for $n=2\cdot7$ or $2\cdot3\cdot7$, respectively, and so $L=\mz_2$ or $\mz_6$. Write $L=\langle x\rangle$. For $n=2\cdot 7$,  we have $\Omega_i=\{\delta_{i1}, \delta_{i2}\}$, $1\leq i\leq 7$, and hence $x=(\delta_{11}, \delta_{12})\cdots (\delta_{71}, \delta_{72})$, an odd permutation on $\Omega$. This is impossible because $T$ is the alternating group $A_{14}$ on $\Omega=\cup_{i=1}^7\Omega_i$. For $n=2\cdot3\cdot7$, we have $\Omega_i=\{\delta_{i1}, \cdots,\delta_{i6}\}$, $1\leq i\leq 7$, and hence we may assume $x=(\delta_{11}, \cdots,\delta_{16})\cdots (\delta_{71},\cdots, \delta_{76})$. Again $x$ is an odd permutation, and this is impossible because $T$ is the alternating group $A_{42}$  on $\Omega$.
\hfill\qed

\medskip
\f {\bf Acknowledgements:} This work was supported by the National Natural Science Foundation of China (11571035, 11731002) and by the 111 Project of China (B16002).

\end{document}